\documentclass{article}
\usepackage{graphicx} 
\usepackage{amsmath, amssymb}

\usepackage[ruled,norelsize]{algorithm2e}
\usepackage{hyperref}
\usepackage{subcaption}
\usepackage[percent]{overpic}
\usepackage{multirow}
\usepackage[top=1in, bottom=1in, left=1in, right=1in]{geometry}

\title{Long Sequence Decoder Network for Mobile Sensing}

\author{Jiazhong Mei$^*$ and J. Nathan Kutz$^{*,\dag}$\\
{\em $^*$ \small Department of Applied Mathematics, University of Washington, Seattle WA 98195}\\
{\em $^\dag$ \small Department of Electrical and Computer Engineering, University of Washington, Seattle WA 98195}}
\date{}

\begin{document}

\maketitle

\begin{abstract}
The reconstruction and estimation of spatio-temporal patterns poses significant challenges when sensor measurements are limited. 
The use of mobile sensors adds additional complexity due to the change in sensor locations over time. 
In such cases, historical measurement and sensor information are useful for better performance, including models such as Kalman filters, recurrent neural networks (RNNs) or transformer models.  However, many of these approaches often fail to efficiently handle long sequences of data in such scenarios and are sensitive to noise.  
In this paper, we consider a model-free approach using the {\em structured state space sequence} (S4D) model as a deep learning layer in traditional sequence models to learn a better representation of historical sensor data.  Specifically, it is integrated with a shallow decoder network for reconstruction of the high-dimensional state space. 
We also introduce a novel initialization of the S4D model using a Butterworth filter design to reduce noise in the inputs. 
Consequently, we construct a robust S4D (rS4D) model by appending the filtering S4D layer before the original S4D structure.
This robust variant enhances the capability to accurately reconstruct spatio-temporal patterns with noisy mobile sensor measurements in long sequence.  Numerical experiments demonstrate that our model achieves better performance compared with previous approaches. Our results underscore the efficacy of leveraging state space models within the context of spatio-temporal data reconstruction and estimation using limited mobile sensor resources, particularly in terms of long-sequence dependency and robustness to noise. 
\end{abstract}

\section{Introduction}

Sensor technologies are ubiquitous across scientific and engineering domains, revolutionizing the way we collect and analyze data. From static installations for environmental monitoring to the emergence of mobile sensors for applications in domains such as autonomous vehicles and wearable health trackers, these sensors play a pivotal role in modern data-driven systems~\cite{Brunton_Kutz_2022,zeng2014convolutional,ren2004trajectory}.  In many cases, measurements of the full state are impossible, impractical, or not even desired. More commonly, limited sensors are used to infer the full characteristic of the system of interest in high dimension from the measurements they collect.  Thus the fundamental mathematical problem is to approximate the full state space from the limited collected data. We consider the problem of state estimation through time sequence measurements from limited mobile sensors by combining a {\em structured state space sequence} (S4D) model with a decoder network, further leveraging a novel initialization scheme and long temporal sequences to produce a robust model with improved performance in comparison with existing methods.

Mobile sensors are becoming more popular and ubiquitous in many applications, for example human biomechanics motion tracking, ocean dynamics monitoring buoys, drone monitoring, and weather balloons~\cite{alvarez2004evolutionary, rosenberg2020predicting, leonard2007collective}.  The mobility of the sensors provide more flexibility and lower cost compared to installing fixed sensors~\cite{peng2014dynamic, ebers2023leveraging, mei2022mobile}. 
However, unlike stationary sensors, state estimation from mobile sensors brings additional challenges.  Traditional techniques, while effective for static sensor arrays, often fall short when applied to mobile sensors operating in dynamic environments.
These models typically employ linear or non-linear mappings from sensor measurements at the current time step to the full system state.   For instance, leveraging the inherent low-rank features of the system, methodologies such as singular value decomposition (SVD), also referred to as proper orthogonal decomposition (POD), identify dominant modes of the system and construct a linear mapping from measurements to high-dimensional state space~\cite{everson1995karhunen, yaglom1967structure, berkooz1993proper, Brunton_Kutz_2022}.  Similarly, dynamic mode decomposition (DMD) extracts linear modes for reconstruction while simultaneously capturing the temporal evolution of these modes in low-rank representation~\cite{Tu2014jcd, kutz2016dynamic, Brunton2022siamreview}.  More complex approaches like shallow decoder networks (SDN) learn nonlinear reconstructions between measurements and high-dimensional state spaces, exhibiting exceptional performance even with a minimal number of sensors~\cite{erichson2020shallow, carter2021data,sahba2022wavefront}.  However, given that the location of each measurement collected by a mobile sensor varies over time, relying solely on such mappings proves inadequate. Considering the impracticality of learning separate models for each sensor location within a high-dimensional state space, there arises a necessity for a generalized model that incorporates sensor location information.  Other approaches, such as the Kalman filter, incorporates historical values alongside current measurements~\cite{stuart2015data}.   
By considering the time history of measurements, additional insights into the system dynamics are gained, enhancing reconstruction performance and robustness to noise. Notably, the Kalman filter naturally accommodates the mobile sensor scenario, where the measurement matrix can vary with the sensor trajectory over time~\cite{mei2022mobile}.  
Despite its adaptability, the Kalman filter is fundamentally a statistical model, necessitating prior knowledge of system dynamics or an approximation, as well as statistical priors regarding noise and disturbances for optimal performance.
Furthermore, the effect of the historical measurements has a compound decay in time depending on the observation noise covariance, lacking the flexibility and ability of memorization in the long run.

The recurrent neural network (RNN)~\cite{rumelhart1986learning} has emerged as a powerful tool for preserving information from past inputs in sequential data commonly seen in a variety of tasks such as speech recognition~\cite{graves2013speech}, machine translation~\cite{sutskever2014sequence}, spatiotemporal predictive learning~\cite{srivastava2015unsupervised, wang2017predrnn}, and much more~\cite{salehinejad2017recent}.
By iteratively applying a series of learnable transformations to input sequences, RNNs adeptly capture temporal dependencies, allowing them to encode and interpret patterns spanning across time. 
A notable approach applying RNN layers to sensing is the {\em SHallow REcurrent Decoder} (SHRED), which has shown promising performance in both stationary and mobile sensor scenarios~\cite{williams2023sensing, ebers2023leveraging}. 
Unlike approaches such as Kalman filter whose performance relies heavily on an approximated model of the system dynamics, SHRED is model-free and directly reconstructs the full system from sensor measurement sequences.
SHRED leverages long short-term memory networks (LSTM)~\cite{hochreiter1997long}, a variant of RNN architecture, in conjunction with a fully-connected, shallow decoder to process time series of sensor measurements for effective reconstruction.
Chen et al.~\cite{chen2020wsn} used a similar deep learning approach combining a recurrent network and a reconstruction network.
Nevertheless, previous research has yet to address several key challenges inherent in mobile sensor reconstruction.

It has been shown that most conventional sequence models such as RNNs and transformers fail to scale to sequences with long time dependencies~\cite{gu2021efficiently, tay2020long}.  
They perform poorly on tasks such as byte-level text classification and retrieval, image classification on sequences of pixels, and finding valid paths connecting two points that are benchmarked by the Long Range Arena~\cite{tay2020long}. 
Long-range dependencies are also very common in a limited mobile sensor reconstruction problem to understand a complex system in a high dimension. 
A mobile sensor would need to take frequent measurements over an extended time to capture the transient and dominant dynamical characteristics of a complex system, which results in a long sequence. 
Therefore, specialized models that address the challenge of long-range dependency should be used.  Second, the model should be robust to sensor failure and disturbances.
Rather than looking at small sensor measurement noise, we are also interested in the case of large measurement error due to failure and disturbance, which could potentially throw off the model and its entire time history memory.
Finally, the study of SHRED on mobile sensors restricted itself to fixed, predetermined sensor trajectories. 
This introduces additional complications in the control of mobile sensors for a  system under the presence of background flows or dynamics.  An excessive amount of energy may be spent to guide the sensors to follow the same trajectory exactly. 
In this work, we aim for a more general model that is independent of the sensor trajectory (unlike mobile SHRED~\cite{ebers2023leveraging}), giving more flexibility and freedom to the sensor control and trajectories.

\begin{figure}[!t]
    \centering
    \begin{minipage}[b]{.5\linewidth}
        \centering
        \begin{overpic}[width=\linewidth]{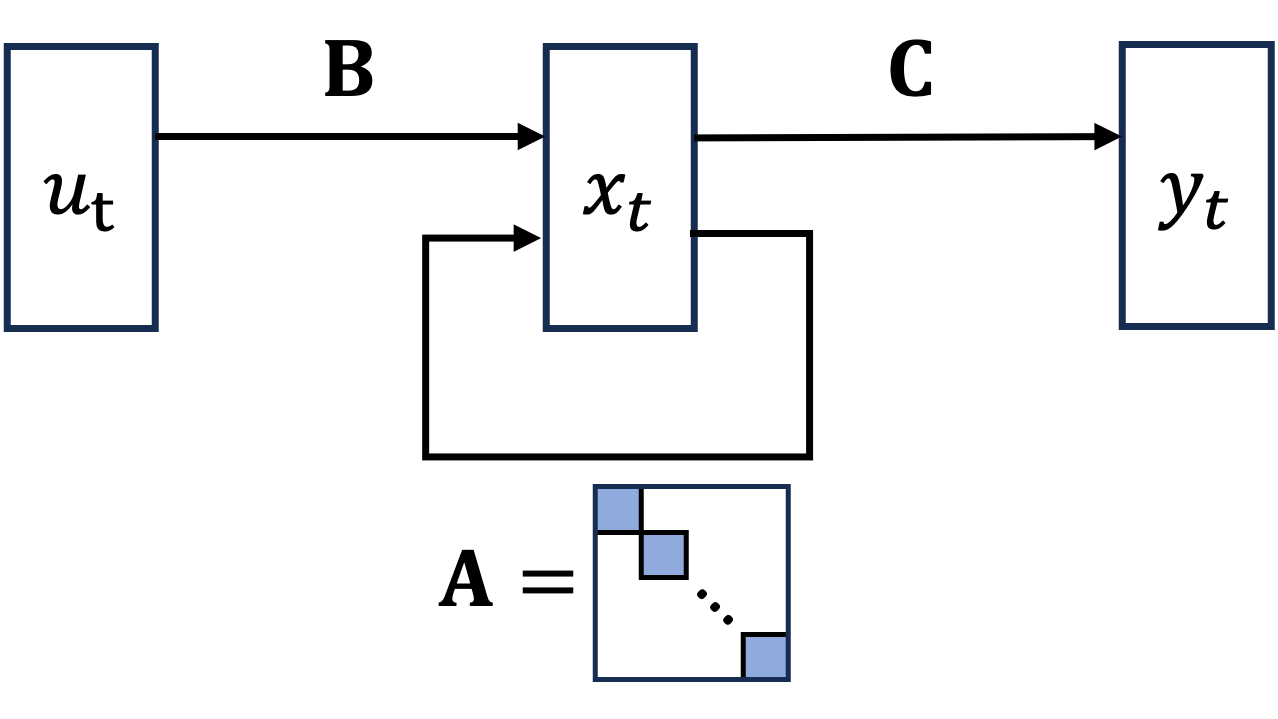}
            \put(0,58){(a)}
        \end{overpic} \\
        \vspace{0.1in}
        \begin{overpic}[width=\linewidth]{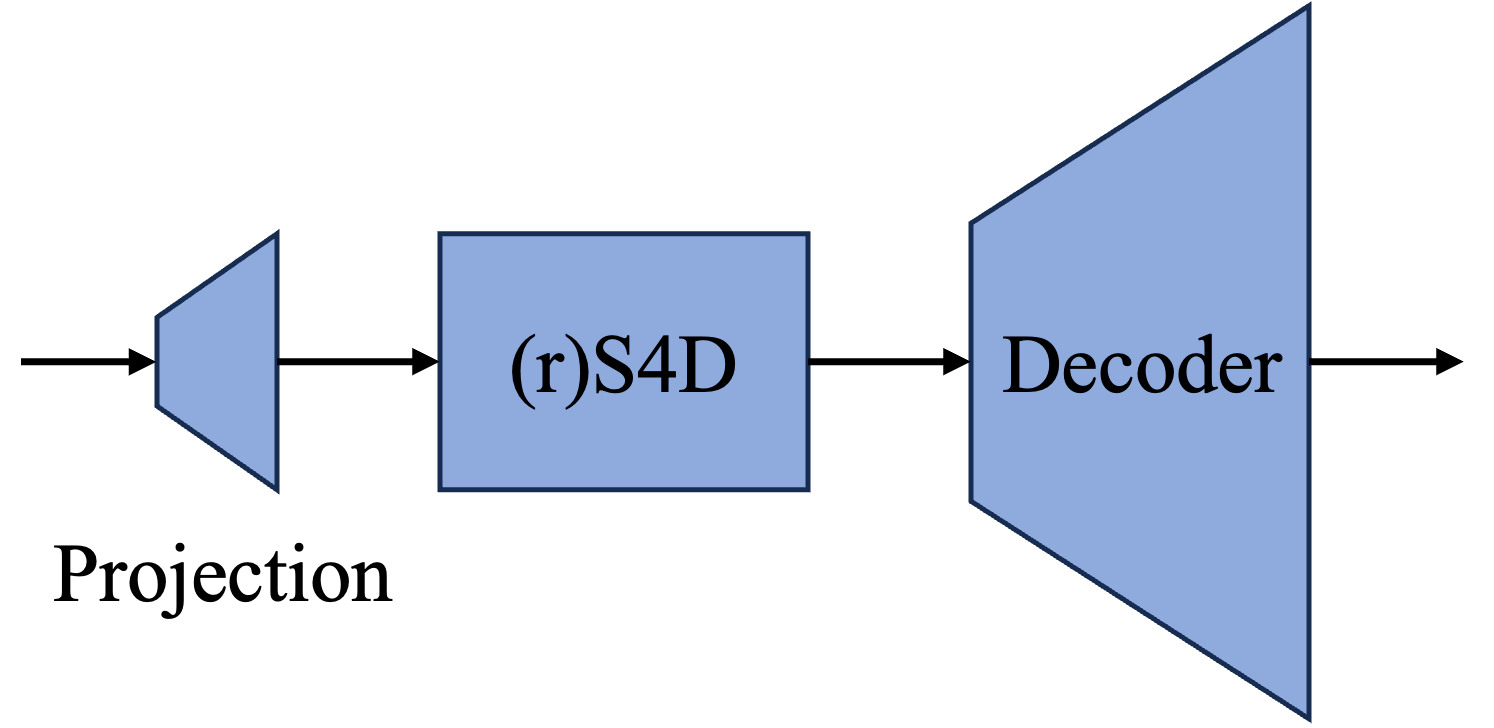}
            \put(0,45){(c)}
        \end{overpic} 
    \end{minipage}%
    \hfill
    \begin{minipage}[t]{.45\linewidth}
        \centering
        \begin{overpic}[width=\linewidth]{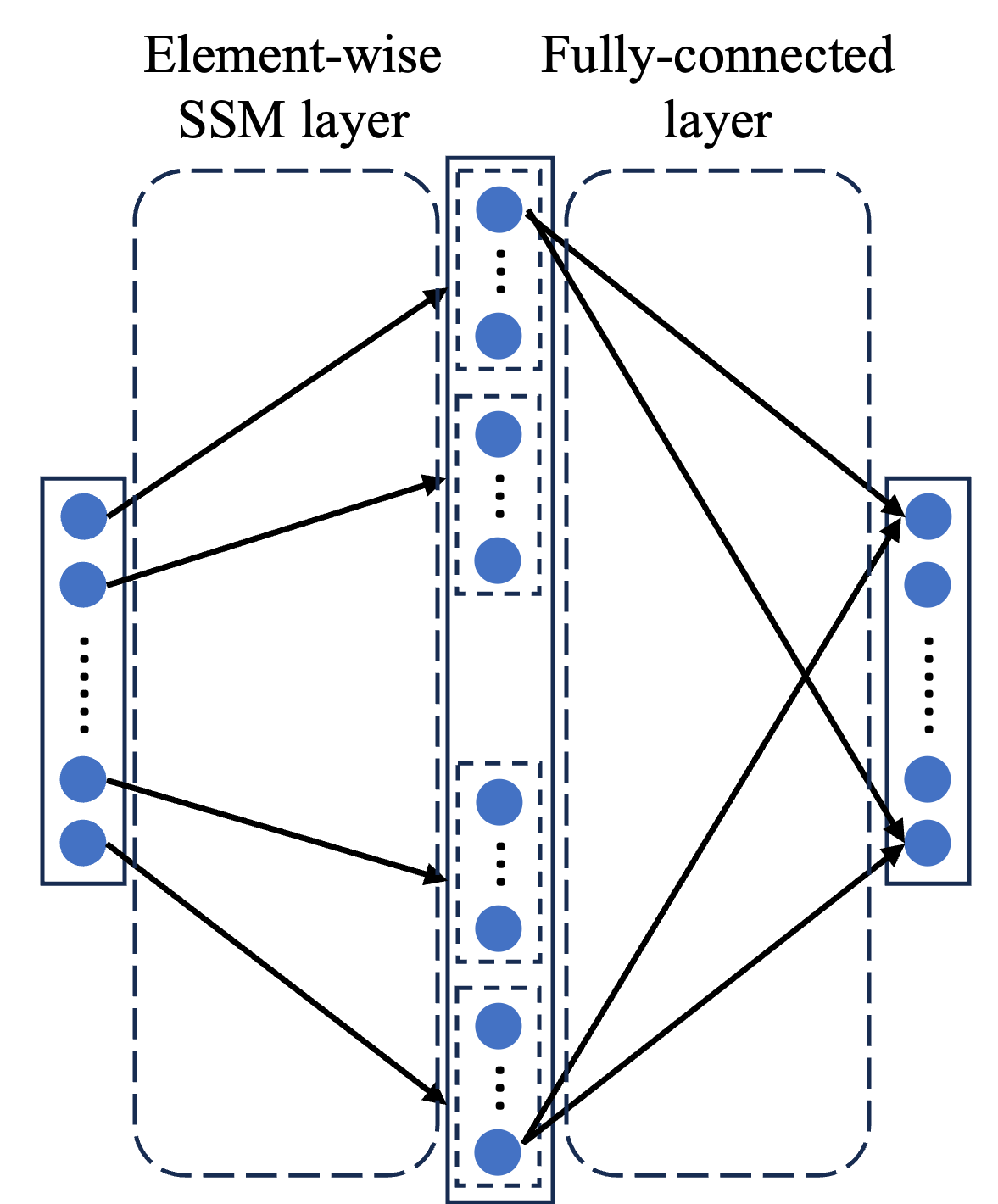}
            \put(-5,101){(b)}
        \end{overpic} 
    \end{minipage}%
    \caption{Model architectures. (a) Standard depiction of {\em state space model} SSM ; (b) {\em structured state space sequence}  (S4D) model block; (c) Integration of {\em shallow recurrent decoder} with S4D architecture to produce the SHRED-(r)S4D model.  The SHRED-(r)S4D model is demonstrated to produce robust and improved performance with arbitrary mobile trajectories.}
    \label{fig:model}
\end{figure}

To address the challenges outlined above, in this paper we propose to employ a state space model in place of the LSTM layers in SHRED as shown in Fig.~\ref{fig:model}.  Specifically, the {\em structured state space sequence}  (S4) model~\cite{gu2021efficiently,gu2022train} is leveraged since it has demonstrated efficacy in handling long-range dependency data, leveraging HiPPO theory for memorization~\cite{gu2020hippo}.
A simplified yet effective variant known as S4D~\cite{gu2022parameterization} utilizes diagonal form approximation to reduce computation complexity and parameterization.
In this study, we introduce a robust variant of the S4D model.
Inspired by the use of HiPPO matrix initialization in S4(D) for long dependency memorization, we propose initializing the SSMs using filtering design to enhance robustness.
Our robust S4D (rS4D) block structure comprises multiple S4D layers, with the first layer initialized using a Butterworth filter for noise filtering and the remaining layers initialized using HiPPO matrix for memorization. 
The rS4D block is seamlessly integrated into the SHRED model, replacing the LSTM block. 
Through numerical experiments, we demonstrate that our model achieves superior performance in long sequence estimation and exhibits reduced sensitivity to measurement noise and disturbances.

\section{Sensing Architecture}

In the following subsections, the mathematical infrastructure is detailed for the proposed architecture of Fig.~\ref{fig:model}.  Specifically, the proposed mobile sensing structure in comprised of various components which when integrated lead to robust and improved performance.

\subsection{State Space Models and S4}

The state space model is defined by the following 1-dimensional input-output, continuous-time, time-invariant system:
\begin{equation}
\begin{split}
    \mathbf{x}'(t) &= \mathbf{A} \mathbf{x}(t) + \mathbf{B} u(t) \\
    y(t) &= \mathbf{C} \mathbf{x}(t)
\end{split}
\end{equation}
where $\mathbf{x}(t) \in \mathbb{R}^N$, $u(t), y(t) \in \mathbb{R}$, $\mathbf{B} \in \mathbb{R}^{N \times 1}$, and $\mathbf{C} \in \mathbb{R}^{1 \times N}$.
It can also be represented in convolution form as:
\begin{equation}
\begin{split}
    K(t) &= \mathbf{C} e^{t\mathbf{A}}  \mathbf{B} \\
    y(t) &= (K \ast u)(t)
\end{split}
\end{equation}
The convolutional form provides computational benefit since it can be converted into a temporal recurrence that is substantially faster for autoregressive applications.  Gu et al.~\cite{gu2021efficiently} showed that the na\"ive state-space model does not work well in practice, possibly due to the fact that the exponential solution to the continuous-time system suffers from vanishing/exploding gradients in a long sequence.
However, the S4 model leverages {\em high-order polynomial projection operators} (HiPPO) theory for parameter initialization and achieves outstanding performance in long-range dependency tasks.

HiPPO theory of continuous-time memorization was first introduced for online function approximation. 
Given a measure that weights the past and some basis functions, HiPPO projects arbitrary functions onto the bases with respect to the measure.
Additionally, the optimal coefficients evolve as a linear ODE with controlling inputs from the target function. 
Therefore, these state coefficients serve as a compressed memorization of the inputs. 
HiPPO matrices refer to a class of transition matrices in the state space models that can memorize the history of input $u(t)$ in the state $\mathbf{x}(t)$.
In particular, a specific HiPPO matrix is defined as follows:
\begin{equation}
    \mathbf{A}_{nk} = -
    \begin{cases}
        (2n+1)^{1/2} (2k+1)^{1/2} & n > k \\
        n+1 & n = k \\
        0 & n < k
    \end{cases}
\end{equation}
The model can be transformed to diagonal plus low-rank (DPLR) form for efficient computation of the convolution.

The structure of the S4 block is then set up as follows. 
A separate state-space model with 1D input and output is considered for each feature element in the multi-dimensional feature input. 
A DPLR dynamics is used, so $\mathbf{A} = \text{diag}(\mathbf{a}) - \mathbf{p}\mathbf{p}^*, \mathbf{a},\mathbf{p} \in \mathbb{R}^N$. 
Along with $\mathbf{B}, \mathbf{C} \in \mathbb{R}^N$ , they make up of the parameters in a S4 block. 
The low rank component can be expressed more generally as the outer product of two separate vectors, but it suffers from numerical instability~\cite{goel2022s}.
Additionally, the output can have multiple values in the form of channels, where $\mathbf{y}(t) \in \mathbb{R}^C$ and $\mathbf{C} \in \mathbb{R}^{C \times N}$. 
The parameters are initialized using the HiPPO matrix. 
Then, the continuous SSM is discretized by a step size $\triangle$ and the outputs are efficiently computing using convolution. 
Suppose we have feature inputs of size $H$. 
S4 handles multiple features by simply defining $H$ independent copies of the state-space model, and then mixing the $CH$ outputs with a position-wise linear layer.
The total number of parameters in a S4 layer is $O(CHN) + O(CH^2)$. 

\subsection{S4D}

Instead of using a DPLR matrix, Gu et al.~\cite{gu2022parameterization} utilize diagonal matrices for further improved efficiency and comparable performance. 
The initialization can be the diagonal component of HiPPO matrix decomposition, or other forms of approximation of the HiPPO matrix.
We can write the kernel as a Vandermonde matrix-vector multiplication, whose discrete convolution kernel depends only on the element-wise product $\mathbf{B} \circ \mathbf{C}_i$ for each channel $i$ of $\mathbf{C}$. 
Therefore, we can train just on $\mathbf{C}$ while keeping $\mathbf{B}=1$ constant.
However, it is shown from experiments that training $\mathbf{B}$ and $\mathbf{C}$ independently gives minor but consistent improvement in performance. 

There are many initializations of the S4D dynamics to approximate tge HiPPO matrix detailed in \cite{gu2022parameterization}.
For example, S4D-LegS directly takes the diagonal values from the HiPPO-LegS matrix; S4D-Inv and S4D-Lin simplify and approximate HiPPO-LegS and HiPPO-FouT matrices with the diagonal values defined as follows:
\begin{equation}
    \textbf{(S4D-Inv)} \quad a_n = -\frac{1}{2} + i\frac{N}{\pi}(\frac{N}{2n+1}-1)
\end{equation}
\begin{equation}
    \textbf{(S4D-Lin)} \quad a_n = -\frac{1}{2} + i\pi n
\end{equation}
In this paper, we use S4D-Lin initialization as it has a simpler form and has shown to be slightly better empirically in the original work.

\subsection{Robust S4D and S4D-BW layer}

An approach to control system robustness and sensitivity to noise is through filtering.
Low pass filtering is a method common is signal processing to remove signals with frequencies higher than a cutoff frequency. 
The Butterworth filter is a type of low pass filter designed to have a frequency response that is as flat as possible in the passband. 
The transfer function of a $N$th-order Butterworth low-pass filter is given by 
\begin{equation}
    G_{BW_N}(s) = \prod_{n=1}^N \frac{\omega_c}{s - \omega_c e^{\frac{i(2n+N-1)\pi}{2N}}},
\end{equation}
where $\omega_c$ is the cutoff frequency.
The poles lie on a circle of radius $\omega_c$ at equally-spaced points, $s_n = \omega_c e^{\frac{i(2n+N-1)\pi}{2N}}$.

S4 models uses HiPPO theory and matrices to initialize the dynamics of the SSMs for memorization.
Then, S4D approximates the dynamics in a diagonal form to promote more efficient computation.
Similarly, if the objective of the SSM changes to act as a low-pass filter, we can set up the initialization differently using the transfer function of the filter.  We consider the general Butterworth filter for low-pass filtering.
The transfer function of the diagonal dynamics in S4D is given by 
\begin{equation}
    G(s) = \mathbf{C} (s\mathbb{I} - \mathbf{A})^{-1} \mathbf{B} = \sum_{n=1}^N \frac{c_nb_n}{s-a_n}.
\end{equation}
It has poles at $s_n = a_n$.
Therefore, we can set the diagonal values to be the poles of the $N$th-order Butterworth low-pass filter.
The order of the filter is represented by the state size of the dynamics.
The cutoff frequency $\omega_c$ is controlled by training the discrete step size $\triangle$ in S4D model. 
We call this S4D-BW. 
\begin{equation}
    \textbf{(S4D-BW)} \quad a_n = e^{\frac{i(2n+N-1)\pi}{2N}}
\end{equation}
S4D-BW can be introduced in front of the regular S4D-Lin layers to filter out high-frequency noise in the inputs before memorization. 
More generally, this opens up possibilities of initialization to the S4D model for filtering purpose, such as high-pass Butterworth filter, and other types of filters. 
The general combined structure of filtering S4D layers and HiPPO S4D layers is called robust S4D (rS4D).

\subsection{SHRED-rS4D for Mobile Sensors}

We append the rS4D block with a shallow decoder network, similar to the SHRED model. 
To differentiate these models, we note the original SHRED model as SHRED-LSTM, a shallow decoder network with a S4D block as SHRED-S4D, and a shallow decoder network with a rS4D block as SHRED-rS4D.
The inputs contain sensor measurements as well as sensor locations. 
The outputs are the full high-dimensional state of the system.
We consider the reconstruction loss in terms of the mean squared error (MSE) loss over the last $t$ time steps.
This is to promote continuous reconstruction performance after some warmup time.
Our code is available at \url{https://t.ly/KGwdA}.

\section{Experimental results}

To demonstrate the SHRED-rS4D method, a number of challenging example problems are considered.  This includes both computational examples of complex spatio-temporal systems (double gyre dynamics, 2D Kolmogorov flow, 2D detonation waves) as well as real data sources (sea-surface temperature data).   The results are compared across architectures.




\subsection{Double Gyre}

A double-gyre flow is a flow pattern that is often seen in many geophysical flows and well studied for its coherent structures~\cite{michini2014robotic,nadiga2001global}. 
We define it using the following stream function
\begin{equation}
\begin{split}
    \psi(x, y, t) = A \sin (\pi f(x,t)) \sin(\pi y), \\ f(x, t) = \epsilon \sin(\omega t)x^2 + x - 2\epsilon \sin(\omega t)x,
\end{split}
\end{equation} 
on a closed and bounded domain $[0,2] \times [0,1]$.
We take the parameters with values $A = 0.5, \omega = 2\pi, \epsilon=0.25$ such that the flow has a period $1$ and a max velocity of $\pi A \approx 1.57$.
We model the vorticity (curl of velocity field) on a $201 \times 101$ discretized grid with a step size of 0.01. 

\begin{figure}[!t]
    \centering
    \includegraphics[width=.7\textwidth]{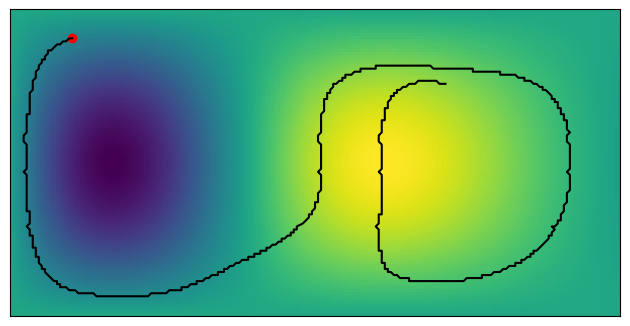}
    \includegraphics[width=.7\textwidth]{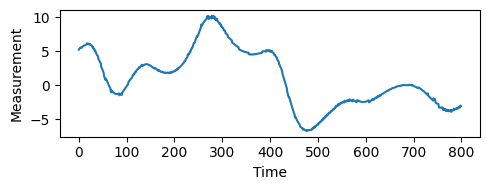}
    \includegraphics[width=.7\textwidth]{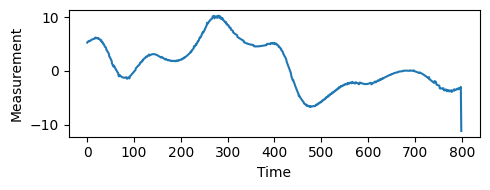}
    \includegraphics[width=.7\textwidth]{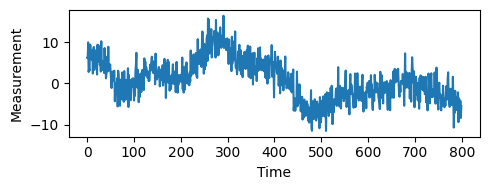}
    \caption{An example of sensor trajectories and measurement inputs collected along the trajectory in the double-gyre system. The line plots are noise-free, disturbed at final step, and noisy measurements.}
    \label{fig:dg_eg}
\end{figure}

We consider one passive mobile sensor floating with the background flow of the system for a total time of $T=4$ covering 4 system periods. 
The flow is given by the velocity field
\begin{equation}
\begin{split}
    \mathbf{v}(x,y,t) &= \begin{bmatrix} -\frac{\partial \psi}{\partial y} \\ \frac{\partial \psi}{\partial x}  \end{bmatrix} \\
    &= \begin{bmatrix} -\pi A \sin(\pi f(x,t)) \cos(\pi y) \\ -\pi A \cos(\pi f(x,t)) \sin(\pi y) \frac{df}{dx} \end{bmatrix} 
\end{split}
\end{equation}
The trajectory of the sensor is varied by initial sensor location and time, both of which are generated randomly.
An example is shown in Figure~\ref{fig:dg_eg}.
Vorticity measurements are collected with a discrete time step of $0.005$ to obtain long sequence dependency.
The vorticity values are standardized and the sensor locations are normalized between 0 and 1.
We generate 2048, 512, 512 random samples for training, validation, and testing repectively.

The models are set as follows. 
We set the main structure of the recurrent component for memorization to be consistent across models for comparison.
That is, the LSTM block in SHRED-LSTM and the S4D-Lin block in SHRED-(r)S4D each contains 2 hidden layers with a hidden dimension of 64. 
The feed-forward decoder component contains 2 layers with a hidden dimension of 128.
MSE loss is computed over the second half of the trajectory, with the belief that enough spatio-temporal information is collected by sensor for reconstruction in the first half. 

\begin{figure}[!t]
    \centering
    \begin{overpic}[height=2in]{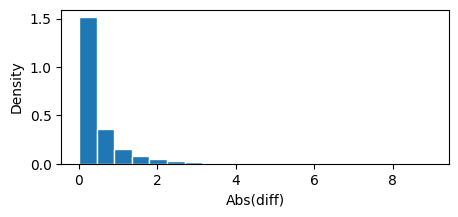}
        \put(-5,42){(a)}
        \put(55,40){$\mu = 0.4827$, $\sigma = 0.6855$}
    \end{overpic}\\
    \begin{overpic}[height=2in]{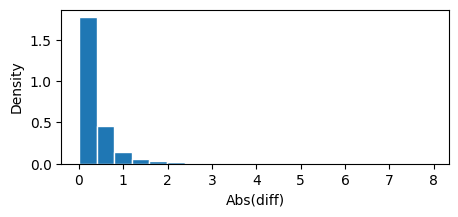}
        \put(-5,42){(b)}
        \put(55,40){$\mu = 0.3704$, $\sigma = 0.5126$}
    \end{overpic}\\
    \begin{overpic}[height=2in]{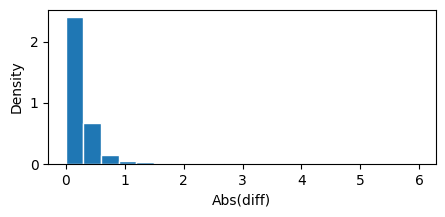}
        \put(-5,42){(c)}
        \put(55,40){$\mu = 0.2565$, $\sigma = 0.3556$}
    \end{overpic}
    \caption{Distribution of the absolute estimation difference at the disturbed time step. (a) SHRED-LSTM; (b) SHRED-S4D; (c) SHRED-rS4D.}
    \label{fig:dg_diff}
\end{figure}

First, we compare the reconstruction performances of the models in a noise-free setting. 
We evaluate the model's performance on normal measurement samples as well as its immediate response to large disturbances, where the sensor measurement at the last time step is severely corrupted. 
The MSE loss is presented in Table~\ref{tbl:dg_loss}.
Notably, SHRED-S4D and SHRED-rS4D demonstrate superior performance over SHRED-LSTM.
Although SHRED-S4D significantly reduces loss in the noise-free test set, it demonstrates high sensitivity to measurement disturbance with the loss more than doubled. 
Conversely, SHRED-rS4D enhances performance in both noise-free and disturbed test sets compared to the benchmarks. 
To delve deeper, we analyze the distribution of the absolute estimation difference at grid points across the disturbed test sets, as depicted in Figure~\ref{fig:dg_diff}. 
It is evident that SHRED-rS4D yields estimations closer to the true system state at the time step when the measurement is corrupted, underscoring its robustness in adverse conditions.

\begin{figure}[!t]
    \centering
    \includegraphics[width=.45\textwidth]{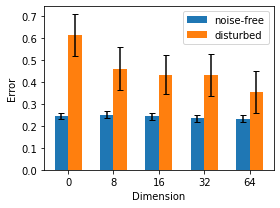}
    \includegraphics[width=.45\textwidth]{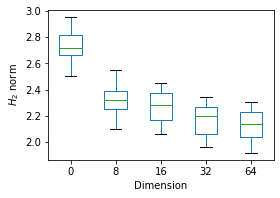} 
    \caption{RMSE error and average $H_2$ norm in the HiPPO layers against the dimension of S4D-BW filtering layer in SHRED-rS4D. 0 dimension refers to the SHRED-S4D model with no filtering.}
    \label{fig:dg_rs4d}
\end{figure}

\begin{figure}[!t]
    \centering
    \begin{overpic}[width=.45\textwidth]{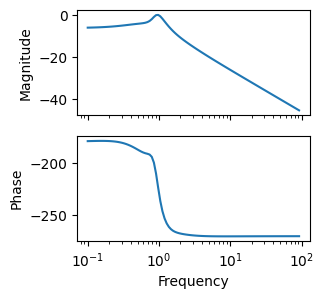}
        \put(-5,85){(a)}
    \end{overpic}
    \hspace{0.2in}
    \begin{overpic}[width=.45\textwidth]{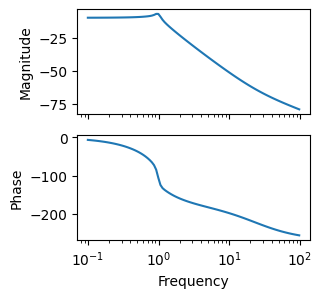}
        \put(-5,85){(b)}
    \end{overpic}
    \vfill
    \begin{overpic}[width=.45\textwidth]{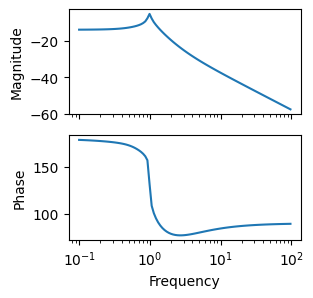}
        \put(-5,90){(c)}
    \end{overpic}
    \hspace{0.2in}
    \begin{overpic}[width=.45\textwidth]{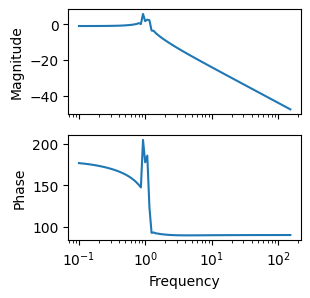}
        \put(-5,90){(d)}
    \end{overpic}\\
    \caption{Bode plots of SSMs in trained S4D-BW filtering layer of dimensions: (a) 8; (b) 16; (c) 32; (d) 64.}
    \label{fig:dg_bode}
\end{figure}

We vary the dimension of the state in the S4D-BW layer to explore potential trade-offs in performance.
Theoretically, while a high-order Butterworth filter offers sharper roll-off and better flatness in the passband, it often suffers from filter instability, such as overshoots and ringing, in step response. 
Thus, we aim to investigate whether a similar trade-off exists in the performance of rS4D and the S4D-BW layer within. 
Our findings, illustrated in Figure~\ref{fig:dg_rs4d}, reveal an improvement in performance as the dimension increases.
Additionally, we examine the bode plots of SSMs in the S4D-BW layer after training. 
Presented in Figure~\ref{fig:dg_bode}, these plots are examples showcasing the sensitivity in the form of ripples around the cutoff frequency, which increases with dimension growth. 
Consequently, while we observe continued performance enhancement with higher filtering dimensions, we believe that the S4D-BW layer should be carefully monitored to assess sensitivity trade-offs based on the dynamics depicted in the bode plots.

The $H_2$ norm of a system serves as an additional metric for system sensitivity to white noise input, representing the average output gain over all frequencies of the input. 
As depicted in Figure~\ref{fig:dg_rs4d}, the inclusion of the S4D-BW filtering layer exerts a positive influence on the $H_2$ norm of the remaining HiPPO S4D-Lin layers. 
Notably, as the state dimension of S4D-BW increases, we observe a corresponding decrease in the average $H_2$ norm within the S4D-Lin layers.

\begin{table}[!t]
\renewcommand{\arraystretch}{1.3}
\caption{Double Gyre testing RMSE using different models as the recurrent block in SHRED. }
\label{tbl:dg_loss}
\centering
\begin{tabular}{c||c|c||c}
\hline
\bfseries Recurrent Block & \bfseries Noise-free & \bfseries Disturbed & \bfseries Noisy \\
\hline
\bfseries LSTM (SHRED) & 0.6361 & 0.9267 & 0.7181 \\
\bfseries S4D & 0.2451 & 0.6141 & 0.2910 \\
\bfseries rS4D & \bf 0.2319 & \bf 0.3559 & \bf 0.2809 \\
\hline
\end{tabular}
\end{table}

\begin{figure}[!t]
    \centering
    \includegraphics[width=.45\textwidth]{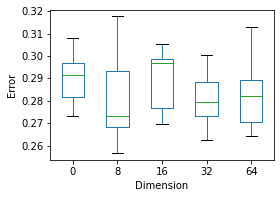}
    \includegraphics[width=.45\textwidth]{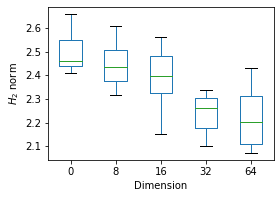} 
    \caption{RMSE error and average $H_2$ norm in the HiPPO layers against the dimension of S4D-BW filtering layer in SHRED-rS4D. 0 dimension refers to the SHRED-S4D model with no filtering.}
    \label{fig:dg_rs4d_noisy}
\end{figure}

Next, we introduce random noise to the sensor measurements throughout the time sequence during training. 
An illustrative example of the sensor measurements along the trajectory is depicted in Figure~\ref{fig:dg_eg}. 
The experimental findings are summarized in Table~\ref{tbl:dg_loss}, revealing that (r)S4D recurrent structures outperform LSTM. 
Interestingly, there is no discernible improvement in performance using rS4D over S4D. 
As depicted in Figure~\ref{fig:dg_rs4d_noisy}, the error remains within a similar range across different dimensions of the filtering layer.
Nonetheless, it is noteworthy that the average $H_2$ norm declines as the dimension increases to 32 and 64. 
This observation suggests that the system's sensitivity improves with additional filtering, rendering the model more resilient to unforeseen disturbances.

\begin{figure}[!t]
    \centering
    \includegraphics[width=.6\textwidth]{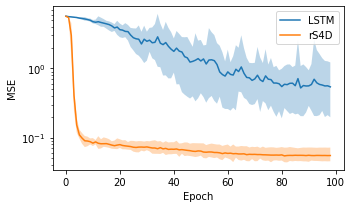}
    \caption{MSE loss vs epochs during training.}
    \label{fig:dg_train}
\end{figure}

Additionally, we note that (r)S4D exhibits significantly better convergence and consistency during training compared to LSTM, as shown in Figure~\ref{fig:dg_train}.
Notably, the MSE loss of SHRED-rS4D rapidly decreases close to the optimal value within a few epochs of training, indicating swift convergence, while SHRED-LSTM learns at a much slower pace. 
Moreover, SHRED-rS4D maintains consistent performance with random model initialization, whereas SHRED-LSTM exhibits higher variance. Consequently, SHRED-rS4D surpasses SHRED-LSTM in terms of model training as well.

\subsection{Sea Surface Temperature}

We study a real-world application regarding the reconstruction of global sea surface temperature from measurements collected by a mobile sensor deployed into the ocean. 
Historically, sensors are deployed in the ocean to provide high-quality observations and validations of the sea surface temperature measured by satellite. 
These sensors can be carried on in situ moorings, drifting buoys, and ships. 
We acquire ocean data from the HYbrid Coordinate Ocean Model (HYCOM) including sea surface temperature, eastward velocity, and northward velocity. 
Daily data is collected on a uniform 1.2-degree lat/lon grid ($134 \times 300$) from 2001 to 2012.

\begin{figure}[!t]
    \centering
    \includegraphics[width=.45\textwidth]{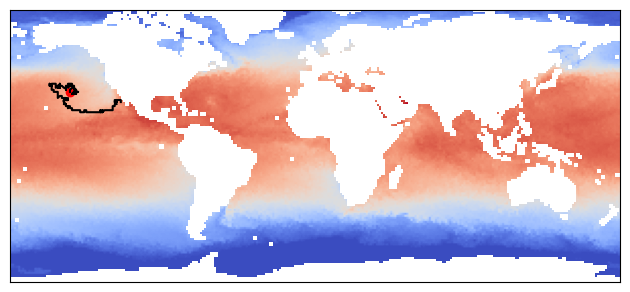}
    \includegraphics[width=.45\textwidth]{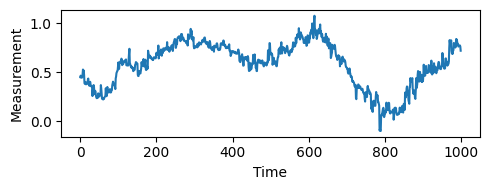}\\
    \includegraphics[width=.45\textwidth]{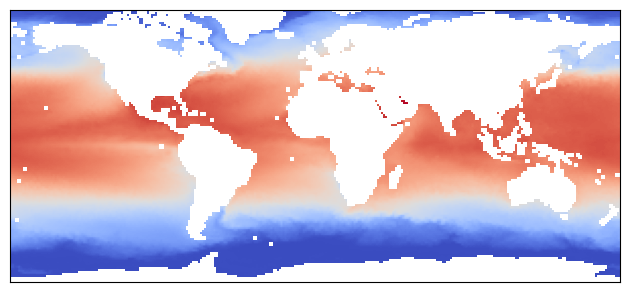}
    \includegraphics[width=.45\textwidth]{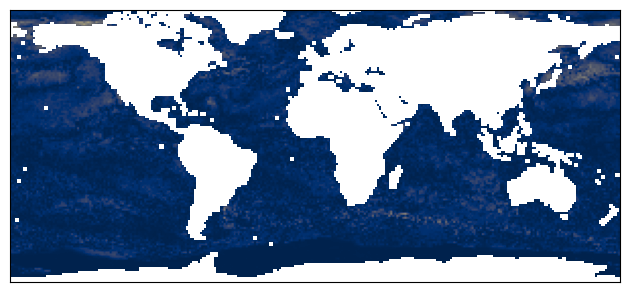}\\
    \includegraphics[width=.4\textwidth]{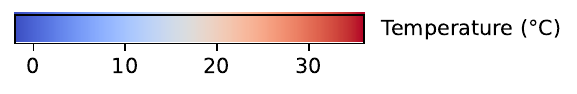}
    \hspace{.1in}
    \includegraphics[width=.4\textwidth]{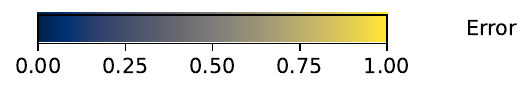}
    \caption{An examples of mobile sensor reconstruction from HYCOM dataset. The top row shows the sensor trajectory and sea surface temperature at the end of the trajectory in the background on the left, and the sensor measurements over time on the right. The bottom row shows the model reconstruction on the left, and the absolute error plotted on the right.}
    \label{fig:hycom_eg}
\end{figure}

The sensor floats passively for a total length of 1000 days with the ocean flow, which is modeled from HYCOM velocity data through continuous interpolation. 
Some examples of the sensor trajectories are shown in Figure~\ref{fig:hycom_eg}.
The daily temperature measurement and the sensor location are used as inputs for the reconstruction of the global sea surface temperature. 
The measurements are standardized and the sensor locations normalized between 0 and 1.
MSE loss is computed for reconstruction on the last day. 
The data is split by time into first 8 years for training and validation and the rest for inference. 
The testing trajectories are such that the reconstruction evaluations are in the inference period. 
10240 trajectories in the training period are generated for training and 2560 trajectories for validation. 
Another 2560 trajectories reaching the inference period are generated for testing. 

The sensor trajectories are generated with random initial time and locations. 
We generated three datasets based on sensor location initialization: random around the globe, random in West Pacific Ocean region, and random in South Atlantic Ocean.
Due to the complex nature of the dynamics of sea surface temperature and the limited coverage of a mobile sensor comparing with the vast ocean space, it is unlikely that a sensor randomly placed on the entire globe guarantees to collect enough information of the system in the given period. 
West Pacific and South Atlantic are two regions that are known to contain comparably richer ocean temperature information than others based on dynamic mode analysis.

The results are summarized in Table~\ref{tbl:hycom_loss}.
We see that SHRED-rS4D has the lowest loss for all regions in validation. 
However, the performance are comparably worse and similar among models in testing, when the estimation time falls out of the training time frame.
This suggests that extrapolation is difficult given the task. 
We also note that for all models, the dataset with trajectories randomly places around the globe has worse performance than those on a more focused regions.
This confirms that the sea surface temperature contains local information and the design of sensor trajectory matters in reconstruction performance.

\begin{table}[!t]
\renewcommand{\arraystretch}{1.3}
\caption{RMSE of HYCOM dataset with random sensor trajectories in different regions.}
\label{tbl:hycom_loss}
\centering
\begin{tabular}{c|c||c|c|c}
\hline
\multirow{2}{*}{\bfseries Type} & \multirow{2}{*}{\bfseries Recurrent Block} & \multicolumn{3}{c}{\bfseries Region} \\ \cline{3-5}
&& \bfseries Global & \bfseries West Pacific & \bfseries South Atlantic\\
\hline
\multirow{3}{*}{\bfseries Val} & \bfseries LSTM (SHRED) & 0.07488 & 0.06176 & 0.06963 \\
& \bfseries S4D & 0.07349 & 0.05939 & 0.06027 \\
& \bfseries rS4D & \bf 0.07241 & \bf 0.05182 & \bf 0.05438 \\
\hline
\multirow{3}{*}{\bfseries Test} & \bfseries LSTM (SHRED) & 0.08276 & 0.08808 & \bf 0.08262 \\
& \bfseries S4D & \bf 0.08225 & \bf 0.08413 & 0.08536 \\
& \bfseries rS4D & 0.08230 & 0.08493 & 0.08605 \\
\hline
\end{tabular}
\end{table}

\subsection{Kolmogorov flow}

We investigate the dynamics of Kolmogorov flow in two dimensions, described by the two-dimensional (2D) Navier–Stokes equations with a sinusoidal body force, across extended periodic domains. The resulting dynamics exhibit the spatio-temporal complexity ideal for challenging data sets. The governing equations are expressed as:
\begin{equation}
    \triangledown \cdot u = 0, \\
    u_t + u \cdot \triangledown u  = -\triangledown p + \nu \triangledown^2 u + f.
\end{equation}
where $u$ represents velocity, $p$ denotes pressure, $\nu$ signifies viscosity, and $f$ stands for the external force term. 
The vorticity is modeled on a discretized $128 \times 128$  grid.

We introduce a single passive mobile sensor, which floats with the background flow of the system for a total duration of 200 seconds, collecting measurements at intervals of 0.1 seconds. 
The system space is wrapped in both directions to accommodate sensor movement. 
Sensor location and start time are initialized randomly. 
We generate 4000, 1000, and 1000 random samples for training, validation, and testing, respectively. 
Testing samples are temporally partitioned from training and validation sets with a later time span.

\begin{figure}[!t]
    \centering
    \includegraphics[width=.7\textwidth]{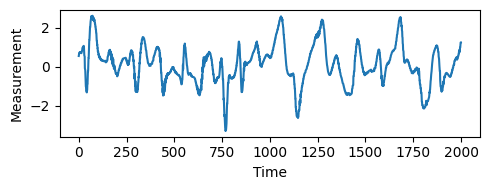}
    \includegraphics[width=.45\textwidth]{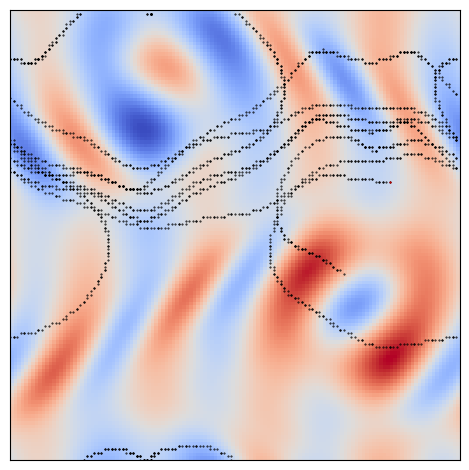} 
    \includegraphics[width=.45\textwidth]{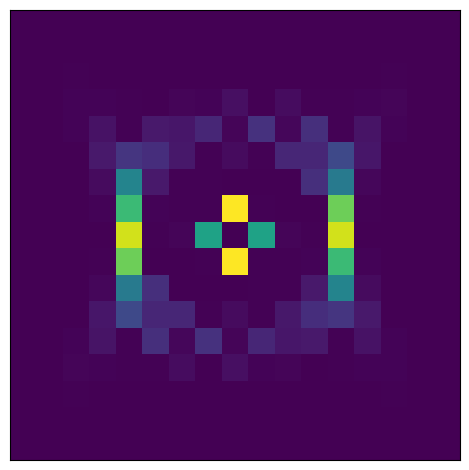}\\
    \includegraphics[width=.45\textwidth]{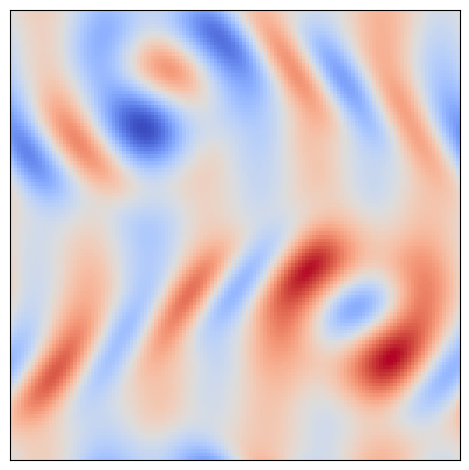} 
    \includegraphics[width=.45\textwidth]{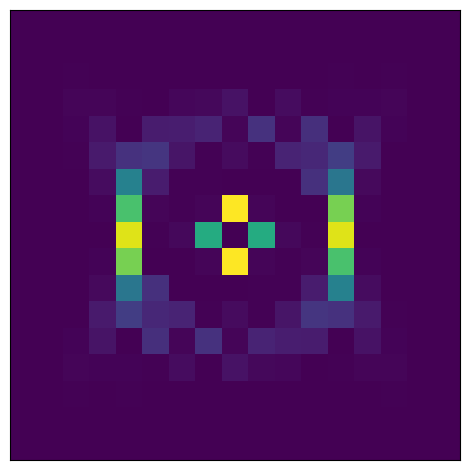}
    \caption{An examples of mobile sensor reconstruction from Kolmogorov dataset. 
    The first row is the sensor measurements along a trajectory.
    The second row shows the true system in physical domain (left) and Fourier domain (right) at the last time step. 
    The sensor trajectory is overlayed on top.
    The thrid row shows the estimation from SHRED-rS4D in physical domain (left) and Fourier domain (right) at the last time step. 
    }
    \label{fig:kol_eg}
\end{figure}

A model comprising 4 recurrent hidden layers and 3 decoder layers is fitted to the data, with mean squared error (MSE) loss computed over the last 100 time steps.
The summarized results in Table~\ref{tbl:summary} and the illustrative estimation example in Figure~\ref{fig:kol_eg} reveal that SHRED-rS4D achieves best performance in validation. 
However, all models exhibit comparable performance in testing. 
Due to the chaotic nature of the system, extrapolation remains challenging.

Interestingly, we observe a pattern in the eigenvalues (diagonal values) of the dynamics matrix in the (r)S4D blocks. 
The trained eigenvalues in the later layers exhibit a much faster decay, except during low oscillation. 
This suggests that most memorization occurs in the earlier layers, with historical patterns of periodicity becoming less useful in later stages for system reconstruction. 
This observation aligns with the non-periodic behavior of the Kolmogorov flow and elucidates the poor performance in extrapolation.

\begin{figure}[!t]
    \centering
    \begin{overpic}[width=.45\textwidth]{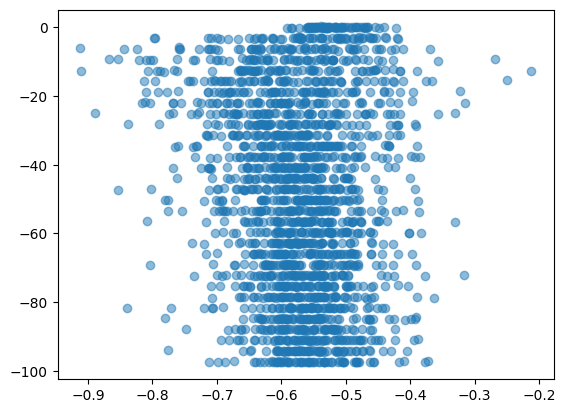}
        \put(10,73){Layer 1}
    \end{overpic} 
    \begin{overpic}[width=.45\textwidth]{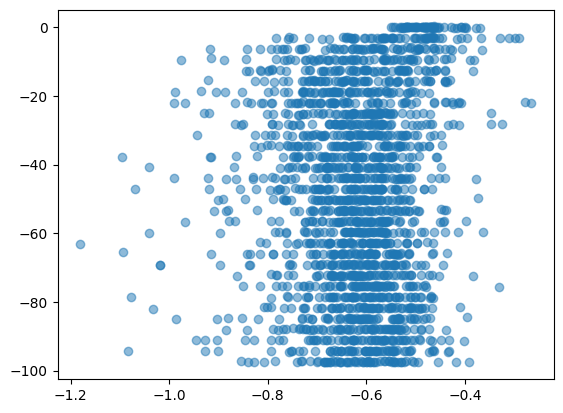}
        \put(10,73){Layer 2}
    \end{overpic} \\
    \vspace{.1in}
    \begin{overpic}[width=.45\textwidth]{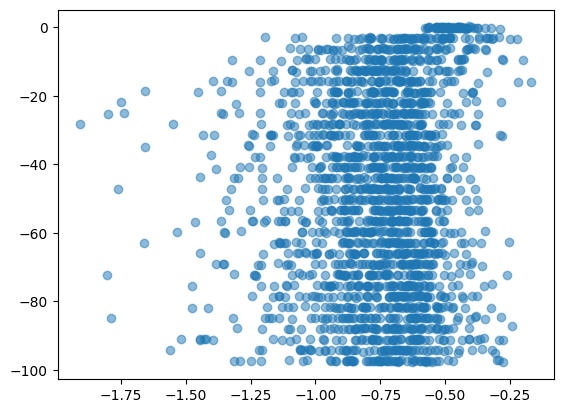}
        \put(10,73){Layer 3}
    \end{overpic} 
    \begin{overpic}[width=.45\textwidth]{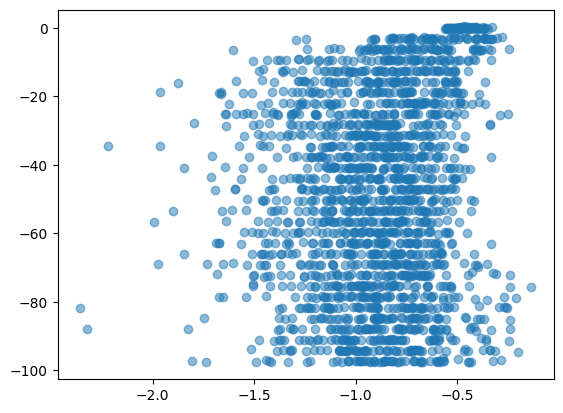}
        \put(10,73){Layer 4}
    \end{overpic} \\
    \caption{Eigenvalues of the dynamics of S4D-Lin layers in SHRED-rS4D.
    }
    \label{fig:kol_param}
\end{figure}

\subsection{Detonation Wave}

We consider the evolution of detonation waves and simulate a variety of explosion scenarios by varying the parameter $\rho_{0_{TNT}}$, which indicates the density of the TNT in the detonation and thus affecting the strength of the explosion:
\begin{align*}
    \rho_{0_{TNT}} = [1000, 1250, 1500, 1650, 1700, 1750, 2000]
\end{align*}
The simulation videos run from $t=0$ to $t=0.001$ with a time step $dt=10^{-6}$. 
We consider the RGB values of the gas concentration as measurements. 

We generate random paths originating close to the center of the explosion and radiating either upward or rightward. 
Each path spans 500 time steps. 
The training set comprises 2000 samples from random trajectories within the first 800 time steps. 
Among these, 80\% are allocated for training, and the remaining 20\% are reserved for validation. 
The testing set consists of 400 samples from random trajectories covering time steps beyond 800.

\begin{figure}[!t]
    \centering
    \includegraphics[width=.7\textwidth]{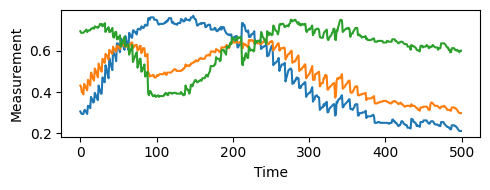}\\
    \includegraphics[width=.45\textwidth]{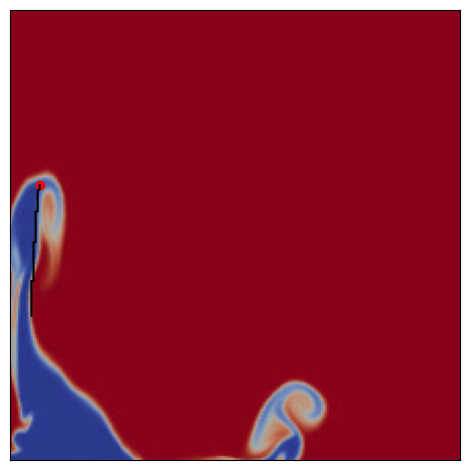}
    \includegraphics[width=.45\textwidth]{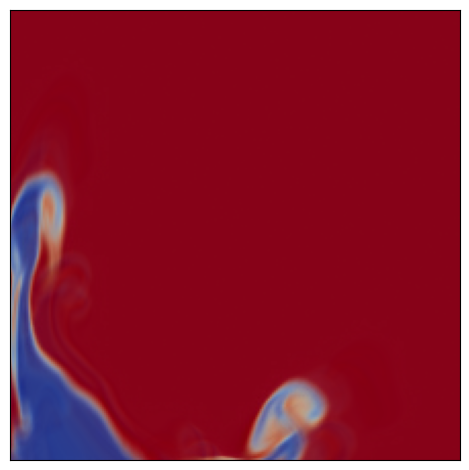} \\
    \includegraphics[width=.45\textwidth]{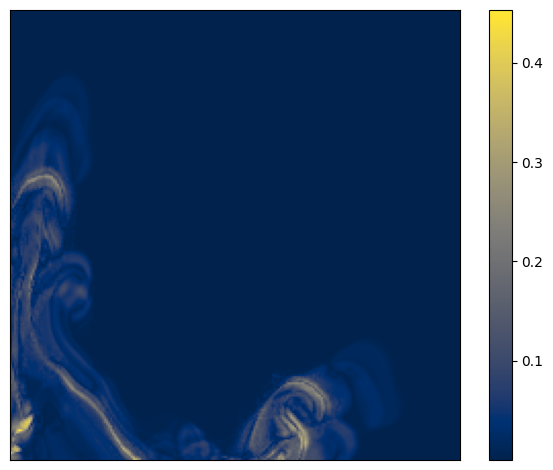}
    \caption{An examples of mobile sensor reconstruction from a small detonation wave dataset. 
    }
    \label{fig:tnt_eg}
\end{figure}

The results are summarized in Table~\ref{tbl:summary}, and an illustrative estimation example is presented in Figure~\ref{fig:tnt_eg}. 
Once again, we observe that SHRED-rS4D achieves best performance in validation, whereas the performance is comparable for all models in testing. 
This suggests that extrapolation in complex systems expanding from a point source is exceedingly challenging for these models, particularly when relying on measurements from a single mobile sensor. 
Additionally, the dynamics evolve differently for various density values $\rho_0$, further complicating estimation. 
We posit that constraining the detonation density parameter and increasing the number of mobile sensors could potentially reduce testing error.

\begin{table}[!t]
\renewcommand{\arraystretch}{1.3}
\caption{RMSE of Kolmogorov flow and small detonation wave datasets.}
\label{tbl:summary}
\centering
\begin{tabular}{c||c|c||c|c}
\hline
\multirow{2}{*}{\bfseries Recurrent Block} & \multicolumn{2}{c||}{\bfseries Kolmogorov}  & \multicolumn{2}{c}{\bfseries SmallDet} \\ \cline{2-5}
& \bfseries Validation & \bfseries Test & \bfseries Validation & \bfseries Test \\
\hline
\bfseries LSTM (SHRED) & 0.2451 & 0.7074 & 0.03578 & 0.07374 \\
\bfseries S4D & 0.1587 & \bf 0.6867 & 0.02996 & \bf 0.07115 \\
\bfseries rS4D & \bf 0.1518 & 0.6967 & \bf 0.02849 & 0.07229 \\
\hline
\end{tabular}
\end{table}

\section{Conclusion}

This paper introduces the SHRED-rS4D model for reconstructing high-dimensional complex systems from limited mobile sensor measurements. 
Building upon the original SHRED model with LSTM as the recurrent structure, our model addresses challenges related to long-range dependency and robustness to noise commonly encountered in mobile sensing problems. 
Leveraging the S4D model and initialization based on HiPPO theory for long sequence memorization, we enhance the model's performance. 
Additionally, we introduce S4D-BW initialization using Butterworth filtering design to reduce sensitivity to noise. 
The integrated structure forms the robust S4D model, replacing the LSTM component in the original SHRED. 
Through numerical examples from complex physical systems, we demonstrate that our model achieves superior performance and exhibits improved training efficiency.  We observe the diminishing effect of the S4D-BW layer as its dimension increases. 
However, extrapolation tasks remain challenging for our model, particularly on non-periodic complex systems. 
We suggest that incorporating more information about the system, such as the number of sensors, length of historical sensor measurements, and sensor coverage, has the potential to reduce extrapolation error.

In this paper, we focus on Butterworth filtering initialization of the S4D model. 
However, future research can explore alternative initialization methods, such as Chebyshev or elliptic filters, to address noise sensitivity more comprehensively. 
Additionally, investigating different objectives for initializing the S4D dynamics could provide insights into addressing other challenges inherent in the reconstruction tasks, which may enhance the model's ability to generalize and improve performance in extrapolation tasks. 

More broadly, the results show the incredible promise of future mobile sensing platforms.  With the ever increasing usage of drone sensing platforms, the SHRED-rS4D architecture allows for a flexible and robust framework that does not necessitate a prescribed trajectory.  As such, it gives a more flexible framework that is not easily compromised by noise or missing/corrupt data.   Moreover, the model is lightweight in comparison with most neural network architectures, thus requiring only modest training data for achieving its superior performance.

\section*{Acknowledgements}
The authors acknowledge support from the Air Force Office of Scientific Research (FA9550-19-1-0386).
The authors further acknowledge support from the National Science Foundation AI Institute in Dynamic Systems
(grant number 2112085).

\bibliographystyle{plain}
\bibliography{bibliography}

\begin{thebibliography}{10}

\bibitem{alvarez2004evolutionary}
Alberto Alvarez, Andrea Caiti, and Reiner Onken.
\newblock Evolutionary path planning for autonomous underwater vehicles in a variable ocean.
\newblock {\em IEEE Journal of Oceanic Engineering}, 29(2):418--429, 2004.

\bibitem{berkooz1993proper}
Gal Berkooz, Philip Holmes, and John~L Lumley.
\newblock The proper orthogonal decomposition in the analysis of turbulent flows.
\newblock {\em Annual review of fluid mechanics}, 25(1):539--575, 1993.

\bibitem{Brunton2022siamreview}
Steven~L Brunton, Marko Budi{\v{s}}i{\'c}, Eurika Kaiser, and J~Nathan Kutz.
\newblock Modern {K}oopman theory for dynamical systems.
\newblock {\em SIAM Review}, 64(2):229--340, 2022.

\bibitem{Brunton_Kutz_2022}
Steven~L. Brunton and J.~Nathan Kutz.
\newblock {\em Data-Driven Science and Engineering: Machine Learning, Dynamical Systems, and Control}.
\newblock Cambridge University Press, 2 edition, 2022.

\bibitem{carter2021data}
Douglas~W Carter, Francis De~Voogt, Renan Soares, and Bharathram Ganapathisubramani.
\newblock Data-driven sparse reconstruction of flow over a stalled aerofoil using experimental data.
\newblock {\em Data-Centric Engineering}, 2:e5, 2021.

\bibitem{chen2020wsn}
Jiahong Chen, Teng Li, Jing Wang, and Clarence~W de~Silva.
\newblock Wsn sampling optimization for signal reconstruction using spatiotemporal autoencoder.
\newblock {\em IEEE Sensors Journal}, 20(23):14290--14301, 2020.

\bibitem{ebers2023leveraging}
Megan~R Ebers, Jan~P Williams, Katherine~M Steele, and J~Nathan Kutz.
\newblock Leveraging arbitrary mobile sensor trajectories with shallow recurrent decoder networks for full-state reconstruction.
\newblock {\em arXiv preprint arXiv:2307.11793}, 2023.

\bibitem{erichson2020shallow}
N~Benjamin Erichson, Lionel Mathelin, Zhewei Yao, Steven~L Brunton, Michael~W Mahoney, and J~Nathan Kutz.
\newblock Shallow neural networks for fluid flow reconstruction with limited sensors.
\newblock {\em Proceedings of the Royal Society A}, 476(2238):20200097, 2020.

\bibitem{everson1995karhunen}
Richard Everson and Lawrence Sirovich.
\newblock Karhunen--loeve procedure for gappy data.
\newblock {\em JOSA A}, 12(8):1657--1664, 1995.

\bibitem{goel2022s}
Karan Goel, Albert Gu, Chris Donahue, and Christopher R{\'e}.
\newblock It’s raw! audio generation with state-space models.
\newblock In {\em International Conference on Machine Learning}, pages 7616--7633. PMLR, 2022.

\bibitem{graves2013speech}
Alex Graves, Abdel-rahman Mohamed, and Geoffrey Hinton.
\newblock Speech recognition with deep recurrent neural networks.
\newblock In {\em 2013 IEEE international conference on acoustics, speech and signal processing}, pages 6645--6649. Ieee, 2013.

\bibitem{gu2020hippo}
Albert Gu, Tri Dao, Stefano Ermon, Atri Rudra, and Christopher R{\'e}.
\newblock Hippo: Recurrent memory with optimal polynomial projections.
\newblock {\em Advances in neural information processing systems}, 33:1474--1487, 2020.

\bibitem{gu2022parameterization}
Albert Gu, Karan Goel, Ankit Gupta, and Christopher R{\'e}.
\newblock On the parameterization and initialization of diagonal state space models.
\newblock {\em Advances in Neural Information Processing Systems}, 35:35971--35983, 2022.

\bibitem{gu2021efficiently}
Albert Gu, Karan Goel, and Christopher R{\'e}.
\newblock Efficiently modeling long sequences with structured state spaces.
\newblock {\em arXiv preprint arXiv:2111.00396}, 2021.

\bibitem{gu2022train}
Albert Gu, Isys Johnson, Aman Timalsina, Atri Rudra, and Christopher R{\'e}.
\newblock How to train your hippo: State space models with generalized orthogonal basis projections.
\newblock {\em arXiv preprint arXiv:2206.12037}, 2022.

\bibitem{hochreiter1997long}
Sepp Hochreiter and J{\"u}rgen Schmidhuber.
\newblock Long short-term memory.
\newblock {\em Neural computation}, 9(8):1735--1780, 1997.

\bibitem{kutz2016dynamic}
J~Nathan Kutz, Steven~L Brunton, Bingni~W Brunton, and Joshua~L Proctor.
\newblock {\em Dynamic mode decomposition: data-driven modeling of complex systems}.
\newblock SIAM, 2016.

\bibitem{leonard2007collective}
Naomi~Ehrich Leonard, Derek~A Paley, Francois Lekien, Rodolphe Sepulchre, David~M Fratantoni, and Russ~E Davis.
\newblock Collective motion, sensor networks, and ocean sampling.
\newblock {\em Proceedings of the IEEE}, 95(1):48--74, 2007.

\bibitem{mei2022mobile}
Jiazhong Mei, Steven~L Brunton, and J~Nathan Kutz.
\newblock Mobile sensor path planning for kalman filter spatiotemporal estimation.
\newblock {\em arXiv preprint arXiv:2212.08280}, 2022.

\bibitem{michini2014robotic}
Matthew Michini, M~Ani Hsieh, Eric Forgoston, and Ira~B Schwartz.
\newblock Robotic tracking of coherent structures in flows.
\newblock {\em IEEE Transactions on Robotics}, 30(3):593--603, 2014.

\bibitem{nadiga2001global}
Balasubramanya~T Nadiga and Benjamin~P Luce.
\newblock Global bifurcation of shilnikov type in a double-gyre ocean model.
\newblock {\em Journal of physical oceanography}, 31(9):2669--2690, 2001.

\bibitem{peng2014dynamic}
Liqian Peng, Doug Lipinski, and Kamran Mohseni.
\newblock Dynamic data driven application system for plume estimation using uavs.
\newblock {\em Journal of Intelligent \& Robotic Systems}, 74(1):421--436, 2014.

\bibitem{ren2004trajectory}
Wei Ren and Randy~W Beard.
\newblock Trajectory tracking for unmanned air vehicles with velocity and heading rate constraints.
\newblock {\em IEEE Transactions on control systems technology}, 12(5):706--716, 2004.

\bibitem{rosenberg2020predicting}
Michael~C Rosenberg, Bora~S Banjanin, Samuel~A Burden, and Katherine~M Steele.
\newblock Predicting walking response to ankle exoskeletons using data-driven models.
\newblock {\em Journal of the Royal Society Interface}, 17(171):20200487, 2020.

\bibitem{rumelhart1986learning}
David~E Rumelhart, Geoffrey~E Hinton, and Ronald~J Williams.
\newblock Learning representations by back-propagating errors.
\newblock {\em nature}, 323(6088):533--536, 1986.

\bibitem{sahba2022wavefront}
Shervin Sahba, Christopher~C Wilcox, Austin McDaniel, Benjamin~D Shaffer, Steven~L Brunton, and J~Nathan Kutz.
\newblock Wavefront sensor fusion via shallow decoder neural networks for aero-optical predictive control.
\newblock In {\em Interferometry XXI}, volume 12223, pages 11--17. SPIE, 2022.

\bibitem{salehinejad2017recent}
Hojjat Salehinejad, Sharan Sankar, Joseph Barfett, Errol Colak, and Shahrokh Valaee.
\newblock Recent advances in recurrent neural networks.
\newblock {\em arXiv preprint arXiv:1801.01078}, 2017.

\bibitem{srivastava2015unsupervised}
Nitish Srivastava, Elman Mansimov, and Ruslan Salakhudinov.
\newblock Unsupervised learning of video representations using lstms.
\newblock In {\em International conference on machine learning}, pages 843--852. PMLR, 2015.

\bibitem{stuart2015data}
Andrew Stuart and Kostas Zygalakis.
\newblock Data assimilation: A mathematical introduction.
\newblock Technical report, Oak Ridge National Lab.(ORNL), Oak Ridge, TN (United States), 2015.

\bibitem{sutskever2014sequence}
Ilya Sutskever, Oriol Vinyals, and Quoc~V Le.
\newblock Sequence to sequence learning with neural networks.
\newblock {\em Advances in neural information processing systems}, 27, 2014.

\bibitem{tay2020long}
Yi~Tay, Mostafa Dehghani, Samira Abnar, Yikang Shen, Dara Bahri, Philip Pham, Jinfeng Rao, Liu Yang, Sebastian Ruder, and Donald Metzler.
\newblock Long range arena: A benchmark for efficient transformers.
\newblock {\em arXiv preprint arXiv:2011.04006}, 2020.

\bibitem{Tu2014jcd}
J.~H. Tu, C.~W. Rowley, D.~M. Luchtenburg, S.~L. Brunton, and J.~N. Kutz.
\newblock On dynamic mode decomposition: theory and applications.
\newblock {\em Journal of Computational Dynamics}, 1(2):391--421, 2014.

\bibitem{wang2017predrnn}
Yunbo Wang, Mingsheng Long, Jianmin Wang, Zhifeng Gao, and Philip~S Yu.
\newblock Predrnn: Recurrent neural networks for predictive learning using spatiotemporal lstms.
\newblock {\em Advances in neural information processing systems}, 30, 2017.

\bibitem{williams2023sensing}
Jan~P Williams, Olivia Zahn, and J~Nathan Kutz.
\newblock Sensing with shallow recurrent decoder networks.
\newblock {\em arXiv preprint arXiv:2301.12011}, 2023.

\bibitem{yaglom1967structure}
AM~Yaglom and VA~Tatarski.
\newblock The structure of inhomogeneous turbulence.
\newblock In {\em Atmospheric Turbulence and Radio Wave Propagation}, pages 166--178. Nauka, 1967.

\bibitem{zeng2014convolutional}
Ming Zeng, Le~T Nguyen, Bo~Yu, Ole~J Mengshoel, Jiang Zhu, Pang Wu, and Joy Zhang.
\newblock Convolutional neural networks for human activity recognition using mobile sensors.
\newblock In {\em 6th international conference on mobile computing, applications and services}, pages 197--205. IEEE, 2014.

\end{thebibliography}

\end{document}